\documentclass[12pt,leqno]{article}
\usepackage{amsmath, amsthm, amsfonts, amssymb, color}
\usepackage{mathrsfs}
\theoremstyle{definition}

\newcommand{\scr}[1]{\mathscr #1}
\definecolor{wco}{rgb}{0.5,0.2,0.3}

\numberwithin{equation}{section} \theoremstyle{remark}

\newcommand{\ua}{\uparrow}

\title{{\bf Robin Heat Semigroup and    HWI Inequality
on  Manifolds with Boundary}\footnote{Supported in
 part by NNSFC(10721091) and the 973-Project.}
}
\author{
{\bf Feng-Yu Wang}\\
\footnotesize{School of Mathematical Sci. and Lab. Math. Com. Sys.,
Beijing Normal
University, Beijing 100875, China}\\
\footnotesize{and}\\ \footnotesize{Department of Mathematics,
Swansea University, Singleton Park, SA2 8PP, UK}\\ \footnotesize{Email: wangfy@bnu.edu.cn;
F.Y.Wang@swansea.ac.uk}}
\begin{document}
\def\R{\mathbb R}  \def\ff{\frac} \def\ss{\sqrt} \def\BB{\mathbb
B}
\def\N{\mathbb N} \def\kk{\kappa} \def\m{{\bf m}}
\def\dd{\delta} \def\DD{\Delta} \def\vv{\varepsilon} \def\rr{\rho}
\def\<{\langle} \def\>{\rangle} \def\GG{\Gamma} \def\gg{\gamma}
  \def\nn{\nabla} \def\pp{\partial} \def\tt{\tilde}
\def\d{\text{\rm{d}}} \def\bb{\beta} \def\aa{\alpha} \def\D{\scr D}
\def\EE{\mathbb E} \def\si{\sigma} \def\ess{\text{\rm{ess}}}
\def\beg{\begin} \def\beq{\begin{equation}}  \def\F{\scr F}
\def\Ric{\text{\rm{Ric}}} \def\Hess{\text{\rm{Hess}}}\def\B{\scr B}
\def\e{\text{\rm{e}}} \def\ua{\underline a} \def\OO{\Omega} \def\sE{\scr E}
\def\oo{\omega}     \def\tt{\tilde} \def\Ric{\text{\rm{Ric}}}
\def\cut{\text{\rm{cut}}} \def\P{\mathbb P} \def\ifn{I_n(f^{\bigotimes n})}
\def\C{\scr C}      \def\aaa{\mathbf{r}}     \def\r{r}
\def\gap{\text{\rm{gap}}} \def\prr{\pi_{{\bf m},\varrho}}  \def\r{\mathbf r}
\def\Z{\mathbb Z} \def\vrr{\varrho} \def\ll{\lambda}
\def\L{\scr L}\def\Tt{\tt} \def\TT{\tt}\def\II{\mathbb I}
\def\i{{\rm i}}\def\Sect{{\rm Sect}}\def\E{\scr E}

\maketitle
\begin{abstract}  Let $M$ be a complete connected Riemannian manifold with boundary
$\pp M$, $Q$ a bounded  continuous function on $\pp M$, and $L=
\DD+Z$ for a $C^1$-vector field $Z$ on $M$. By using the reflecting
diffusion process generated by $L$ and its local time on the
boundary, a probabilistic formula is presented for the semigroup
generated by $L$ on $M$ with Robin boundary condition $\<N,\nn
f\>+Qf=0,$ where $N$ is the inward unit normal vector field of $\pp
M$. As an application, the   HWI inequality   is established on
manifolds with (nonconvex) boundary. In order to study this
semigroup, Hsu's gradient estimate and the corresponding Bismut's
derivative  formula are established on a class of noncompact
manifolds with boundary.
\end{abstract} \noindent
 AMS subject Classification:\ 60J60, 58G32.   \\
\noindent
 Keywords:   Gradient estimate, HWI inequality, local time.
 \vskip 2cm

\section{Introduction}
Let $M$ be a $d$-dimensional connected complete Riemannian manifold
with boundary $\pp M$ and $L=\DD+ Z$ for some $C^1$-vector field $Z$
such that

\beq\label{C} \Ric-\nn Z  \ge -K\end{equation} holds on $M$ for some
constant $K\in \R$. This curvature condition is well known by Bakry
and Emery   \cite{BE}.

 Let $X_t$ be the reflecting diffusion process generated by
$L$ on $M$, and let $l_t$ be its local time on the boundary $\pp M.$
Let $\tau$ be the first hitting time of $X_t$ to $\pp M$. It is well
known that the following heat equation can be described by using the
process $X_t$:

\beq\label{H} \pp_t u= Lu,\ \ \ u(0,\cdot)=f,\end{equation} where
$f\in \B_b(M).$ With Dirichlet boundary condition $u|_{\pp M}=0$ the
solution can be formulated as

$$u(t,x)= \EE^x \big[f(X_t)1_{\{t<\tau\}}\big]$$ while under the Neumann
boundary condition $Nu|_{\pp M}=0$ one has

$$u(t,x)= \EE^x f(X_t),$$ where $\EE^x$ is the expectation taking for
the process $X_t$ starting at $x.$ In this paper we shall provide
the corresponding probability formula for the solution under the
Robin boundary condition (cf. \cite[page 102]{O}):

\beq\label{Robin} Nf:= \<N,\nn f\>= -Qf\ \ \text{on}\ \pp
M,\end{equation}where $Q\in C_b(\pp M)$ and $N$ is the inward unit
normal vector field on $\pp M$. It turns out that under a reasonable
assumption the solution to (\ref{H}) under condition (\ref{Robin})
can be formulated by

\beq\label{A0} u(t,x)=P_t^Q f(x):= \EE^x \big\{f(X_t)\e^{\int_0^t
Q(X_s)\d l_s}\big\},\ \ \ t\ge 0, x\in M.\end{equation} As soon as
$P_t^Q$ is well defined, the semigroup property follows immediately
from the Markov property of the reflecting diffusion process $X_t$.
To ensure the boundedness of $P_t^Q$ under the uniform norm, it is
natural to ask the local time $l_t$ to be exponentially integrable.
According to calculations from \cite{W05} (see also the proof of
Lemma \ref{L2.2} below), for this we shall need the following
assumption.

\paragraph{(A)}  \emph{The boundary $\pp M$ has a bounded second fundamental form and a strictly positive
 injectivity radius,   the
sectional curvature of $M$ is bounded above,    and there exists
$r>0$ such that $Z$ is bounded on the $r$-neighborhood of $\pp M$.
}

\ \

Let $\rr_{\pp M}$ be the Riemannian distance to the boundary. Then
the $r$-neighborhood of $\pp M$ is $\pp_r M:= \{x\in M: \rr_{\pp
M}(x)<r\},$ where $\rr_{\pp M}$ is the Riemannian distance to the
boundary $\pp M$.  Next, the injectivity radius $\i_{\pp M}$ of $\pp
M$ is the largest number $r$ such that the exponential map

$$[0,r)\times \pp M\ni (s,x)\mapsto \exp[s N_x]\in \pp_rM$$
is diffeomorphic. In particular, $\rr_{\pp M}$ is smooth on $\pp_r
M$ for $r\le\i_{\pp M}.$

Finally, to state our result,  we introduce the following class of
references functions:

$$ \D_0:=\{f\in C_0^\infty(M):\ Nf+Qf=0\ \text{on}\
\pp M\}.$$

 \beg{thm}\label{T1.1} Assume {\bf (A)} and $(\ref{C})$
hold and let $Q\in C_b(\pp M)$.

$(1)$  $\{P_t^Q\}_{t\ge 0}$ is a positivity-preserving strong Feller
semigroup of bounded linear operators on $\B_b(M)$, whose generator
is $L$ with domain containing all functions $f\in C^2_0(M)$ such
that $Nf+Qf=0$ holds on $\pp M$.

$(2)$ Let $Z=\nn V$ for some $V\in C^2(M)$ with $\mu(\d
x):=\e^{V(x)}\d x$ not necessarily finite. Then $\{P_t^Q\}_{t\ge 0}$
provides  a bounded
 symmetric $C_0$-semigroup on $L^2(\mu).$ If in particular $Q\le 0$, then $P_t^Q$ is sub-Markovian and the associated
 symmetric Dirichlet form  $(\scr E,\D(\scr E))$  is the closure of
 $(\E,\D_0)$ with

$$\scr E(f,g)= \mu(\<f,g\>) -\mu_\pp(Qfg),\ \ \ f,g\in \D_0.$$
  \end{thm}

We note that a solution to the heat equation (\ref{H}) under the
Robin condition (\ref{Robin}) can be represented by (\ref{A0})
provided it is bounded in $x\in M$. Indeed, by the boundary
condition and the It\^o formula, for fixed $t>0$,

$$\d u(t-s, X_s)= \d M_s +Nu(t-s,\cdot)(X_s)\,\d l_s= \d M_s -
(Qu(t-s,\cdot))(X_s)\, \d l_s$$ holds for some local martingale
$M_s$ up to time $t$. So, $s\mapsto u(t-s, X_s)\exp[\int_0^s
Q(X_r)\d l_r]$ is a local martingale as well. Since $u$ is bounded
and $l_t$ is exponential integrable due to Lemma \ref{L2.2} below,
it is indeed a martingale. Thus, (\ref{A0}) holds.

Since the local time $l_t$ is not absolutely continuous in $t$, the
semigroup $P_t^Q$ is essentially different from the well developed
Schr\"odinger semigroup. According to Theorem \ref{T1.1} below
$P_t^Q$ is generated by $L$ under the boundary condition $Nf+Qf=0$,
where $N$ is the inward unit normal vector field on $\pp M$. So, the
formula (\ref{A0}) will be important in the study of this boundary
value problem on $M$. In this paper we shall explain how can one
apply   this semigroup to the study of  HWI inequality on manifolds
with boundary. This inequality links three important quantities
including the entropy, the energy and the Wasserstein distance (or
the optimal transportation cost), and was found in \cite{OV, BGL} on
manifolds without boundary.

To study the HWI inequality,  we consider the symmetric case that
$Z=\nn V$ for some $V\in C^2(M)$ such that $\mu(\d x)=\ \e^{V(x)}\d
x$ is a probability measure on $M$, where $\d x$ is the Riemannian
volume measure on $M$. Let $P_t$ be the semigroup of the reflecting
diffusion process generated by $L$ on $M$, which is  then symmetric
in $L^2(\mu)$. When $\pp M$ is convex (\ref{C}) implies the
following gradient estimate (cf. \cite{Q, W97})

\beq\label{G} |\nn P_t f|\le \e^{Kt}P_t|\nn f|,\ \ \ f\in
C_b^1(M).\end{equation} Combining this estimate and an argument of
\cite{BGL}, we can easily obtain the following HWI inequality:

\beq\label{HWI} \mu(f^2\log f^2)\le 2 \ss{\mu(|\nn f|^2)}\,
W_2(f^2\mu,\mu)+\ff K 2 W_2(f^2\mu, \mu)^2,\ \ \
\mu(f^2)=1,\end{equation} where $W_2$ is the $L^2$-Wasserstein
distance induced by the Riemannian distance function $\rr$ on $M$.
More precisely, for a probability measure $\nu$ on $M$  (note that
we are using $\rr^2$ to replace $\ff 1 2 \rr^2$ in \cite{BGL})

$$W_2(\nu,\mu)^2:= \inf_{\pi\in\scr C(\nu,\mu)} \int_{M\times M} \rr(x,y)^2\pi(\d x,\d y),$$ where
$\scr C(\nu,\mu)$ is the class of all couplings of $\nu$ and $\mu.$

To see that $P_t^Q$ is important in the study of the HWI inequality
on a nonconvex manifold, let us briefly introduce the main idea for
the proof of (\ref{HWI}) on a convex manifold using (\ref{G}).
Firstly, due to Bakry and Emery,
  (\ref{G}) implies the semigroup log-Sobolev inequality

  \beq\label{LS1} P_t( f^2\log f^2)\le (P_t f^2)\log P_t f^2 +\ff{2(\e^{2Kt}-1)}{K} P_t |\nn f|^2.\end{equation}
  Taking integration
  for both sides with respect to $\mu$ we arrive at

  \beq\label{LS2} \mu(f^2\log f^2)\le \ff{2(\e^{2Kt}-1)}{K} \mu(|\nn f|^2) +\mu((P_t f^2)\log P_t f^2)).
  \end{equation}
  On the other hand,  according to \cite[Proof of Lemma 4.2]{BGL} the gradient estimate (\ref{G}) implies
   (again note that the $W_2^2$ here is twice of the one in \cite{BGL})

  \beq\label{W} \mu((P_t f^2)\log P_t f^2)\le \ff{K\e^{2Kt}}{2(\e^{2Kt}-1)} W_2(f^2\mu,\mu)^2.\end{equation}
  Combining this with (\ref{LS2}) and minimizing in $t>0$, one derives  (\ref{HWI}).

  Now, what can we do for the nonconvex setting? According to \cite{Hsu}, in this case the local time and
  the second fundamental form  will be naturally involved in the upper bound of $|\nn P_t f|$.
   Let the second fundamental form be bounded below by $-\si$ for some $\si\ge0,$ i.e.

  \beq\label{F} \II(X,X):= -\<\nn_X N, X\>\ge -\si|X|^2,\ \ \ X\in T\pp M.\end{equation} Recall that
   $N$ is the inward
  unit normal vector field on $\pp M$. According to \cite[Theorem 5.1]{Hsu},
  if $M$ is compact and $V=0$ then (\ref{C}) and (\ref{F}) imply

  \beq\label{G2} |\nn P_t f|(x)\le \e^{Kt} \mathbb E^x \big\{|\nn f|(X_t) \e^{\si l_t}\big\},\ \ \ x\in M, t\ge 0, f\in C_b^1(M),\end{equation}
  where $X_t$ is the reflecting $L$-diffusion process and $l_t$ is its local time on $\pp M$.
In this paper we shall prove (\ref{G2})  for $Z\ne 0$ on noncompact
manifolds under assumption {\bf (A)}, see Proposition \ref{P2.1}
below.

 Since the local time is unbounded, we are not
able to derive from (\ref{G2}) the semigroup log-Sobolev inequality
  like  (\ref{LS1}).
But Theorem \ref{T1.1}(2) enables us to derive a log-Sobolev
inequality of type (\ref{LS2}) using
 (\ref{G2}), from which we can prove the following HWI inequality (\ref{HWI2}).

\beg{thm}\label{T1.2} Let $Z=\nn V$ for some $V\in C^2(M)$ such that
$\mu$ is a probability measure. Assume {\bf (A)} and $(\ref{C})$.
Let $\II\ge -\si$ for some $\si\in \R$. Then

$$\eta_\ll(s):= \sup_{x\in M}\mathbb E^x \e^{\ll l_s}<\infty,\ \ \ s,\ll\ge 0$$ holds, and for any $t>0$,

\beq\label{HWI2} \mu(f^2\log f^2)\le 4\bigg(\int_0^t \e^{2Ks}\eta_{2\si}(s)\d s\bigg)\mu(|\nn f|^2)
+ \ff {W_2(f^2\mu,\mu)^2} {4\int_0^t\e^{-2Ks}\eta_{2\si }(s)^{-1} \d s},\ \ \mu(f^2)=1.\end{equation}
\end{thm}

To derive an explicit HWI inequality, we shall estimate
$\eta_{2\si}$ as in \cite{W05} by using the It\^o formula for
$\phi\circ\rr_{\pp M}(X_t)$ with a specific choice of $\phi$ (see
Lemma \ref{L2.2} below).  From this we obtain the following
consequence of Theorem \ref{T1.1} immediately. Let $\Sect_M$ be the
sectional curvature of $M$, and  let

$$\dd_r(Z):= \sup_{\pp_r M}\<Z, \nn\rr_{\pp M}\>^-,\ \ \ r>0.$$

\beg{cor} \label{C1.2} Let $Z=\nn V$ for some $V\in C^2(M)$ such
that $\mu$ is a probability measure. Assume {\bf (A)} and
$(\ref{C})$. Let $r_0,\si, k,
>0$ be such that $\dd_{r_0}(Z)<\infty, -\si\le\II\le \gg$ and
$\Sect_M\le k$. For any

$$0<r\le \min\bigg\{\i_{\pp M},\ r_0,\ \ff1 {\ss k} \arcsin\bigg(\ff{\ss k}{\ss{k+\gg^2}}\bigg)\bigg\},$$ and for
$K_r:= K + \ff {\si d}{r} +  \si  \dd_r(Z) + 4\si^2,$ the HWI
inequality
$$\mu(f^2\log f^2)\le 2\e^{2\si d/r}\ss{\mu(|\nn f|^2)}\,
W_2(f^2\mu,\mu)+\ff {K_r\e^{2\si d/r}} 2 W_2(f^2\mu, \mu)^2,\ \ \ \mu(f^2)=1$$ holds.\end{cor}

As preparations, in the next section we shall confirm the
exponential integrability of $l_t$ and establish (\ref{G2}) on
noncompact manifolds. The above two theorems are then proved in
Sections 3 and 4 respectively. To prove the strong Feller property
of $P_t^Q$  and for further applications in the literature,  the
Bismut type formula for $P_t$ on manifolds with boundary is
addressed in Appendix at the end of the paper.

\section{Exponential estimate and Hsu's gradient estimate}

As explained in Section 1, to ensure that $P_t^Q$ is well defined,
we first study the exponential integrability of the local time.

\beg{lem} \label{L2.2} Let $r_0>0$ be such that
$\dd_{r_0}(Z)<\infty$ and let $k,\gg$ be in Corollary \ref{C1.2}.
  Then

  $$\sup_{x\in M} \EE^x \e^{\ll l_t} \le \exp\Big[ \ff{\ll d r} 2
  +\Big(\ff{ \ll d} r + \ll \dd_r (Z)+ 2 \ll^2\Big)t\Big],\ \ \
  t\ge 0, \ll\ge 0$$ holds for any

 $$0<r\le \min\bigg\{\i_{\pp M},\ r_0,\ \ff1 {\ss k} \arcsin\bigg(\ff{\ss k}{\ss{k+\gg^2}}\bigg)\bigg\}.$$ \end{lem}

 \beg{proof}Let

 $$h(s)= \cos\big(\ss k\, s\big)- \ff{\gg } {\ss k}\sin\big(\ss k\, s\big),\
 \ \ s\ge 0.$$ Then $h$ is the unique solution to the equation

 $$h'' + k h=0,\ \ \ h(0)=1, h'(0)= -\gg.$$
 By the Laplacian comparison theorem for $\rr_{\pp M}$ (cf.
 \cite[Theorem 0.3]{Kasue} or \cite{W07}),

 $$ \DD \rr_{\pp M} \ge \ff{(d-1)h'}{h}(\rr_{\pp
 M}),\ \ \ \rr_{\pp M}< \i_{\pp M} \land h^{(-1)}(0).$$
 Thus,

 \beq\label{L0} L\rr_{\pp M}\ge \ff{(d-1)h'}{h}(\rr_{\pp M})
 -\dd_r(Z),\ \ \ \rr_{\pp M}\le r.\end{equation} Now, let

 \beg{equation*}\beg{split} &\aa= (1-h(r))^{1-d}\int_0^r
 (h(s)-h(r))^{d-1}\d s,\\
 &\psi(s) =\ff 1 \aa\int_0^s (h(t)-h(r))^{1-d}\d t\int_{t\land
 r}^r(h(u)-h(r))^{d-1}\d u,\ \ s\ge 0.\end{split}\end{equation*} We
 have $\psi(0)=0, 0\le \psi'\le \psi'(0)=1.$ Moreover, as observed in \cite[Proof of
 Theorem 1.1]{W05},

 \beq\label{L1} \aa\ge \ff r d,\ \ \psi(\infty)=\psi(r)\le
 \ff{r^2}{2\aa}\le \ff{dr} 2.\end{equation} Combining this with
 (\ref{L0}) we obtain (note that $\psi'(s)=0$ for $s\ge r$)

 \beq\label{L2} L\psi\circ\rr_{\pp M} =\psi'\circ\rr_{\pp
 M}L\rr_{\pp M}
 +\psi''\circ\rr_{\pp M}\ge -\ff 1 \aa -\dd_r(Z)\ge -\ff d r
 -\dd_r(Z).  \end{equation} On the other hand, since $\psi'(0)=1$, by the It\^o
 formula we have

 \beq\label{CC}\d \psi\circ\rr_{\pp M}(X_t) =\ss 2 \psi'\circ\rr_{\pp M}(X_t) \d
 b_t +L\psi\circ\rr_{\pp M}(X_t)\d t +\d l_t,\end{equation}
 where $b_t$ is the one-dimensional Brownian motion. Then it follows from
 (\ref{L1}) and (\ref{L2}) that (note that $|\psi'|\le 1$)

\beg{equation*}\beg{split} \EE\e^{\ll l_t} &=\EE \exp\bigg[\ll
\psi\circ\rr_{\pp M}(X_t) +\Big(\ff{d\ll} r+\ll \dd_r(Z)\Big)t-\ss 2
\ll \int_0^t\psi'\circ\rr_{\pp M}(X_s)\d b_s\bigg]\\
&\le \exp\Big[\ff 1 2 \ll d r +\Big(\ff{d\ll} r
+\ll\dd_r(Z)\Big)t\Big]
\bigg(\EE\exp\bigg[4\ll^2\int_0^t\big(\psi'\circ\rr_{\pp
M}(X_s)\big)^2\d s\bigg]\bigg)^{1/2}\\
&\le \exp\bigg[\ff 1 2 \ll d r +\Big(\ff {d\ll} r +\ll \dd_r(Z) + 2
\ll^2\Big)t\bigg].\end{split}\end{equation*}
 \end{proof}

 This Lemma ensures  the boundedness of $P_t^Q$ under the
 uniform norm. Next, we intend to prove (\ref{G2}) under assumption {\bf (A)}, which is known by \cite{Hsu} for compact $M$
and $Z=0.$

\beg{prp}\label{P2.1} Assume that {\bf (A)}.  Let $\kk_1, \kk_2\in
C_b(M)$ be such that

\beq\label{2.0} \Ric -\nn Z\ge -\kk_1,\ \ \ \II\ge
-\kk_2\end{equation}  hold on $M$ and $\pp M$ respectively. Then

\beq\label{Hsu} |\nn P_t f|(x) \le \EE^x \bigg\{|\nn f|(X_t)
\exp\bigg[ \int_0^t \kk_1(X_s)\d s+\int_0^t\kk_2(X_s)\d
l_s\bigg]\bigg\}\end{equation}
 holds for all $f\in C_b^1(M), t>0, x\in M.$ \end{prp}

We first provide a simple proof of (\ref{Hsu}) under a further  condition that
 $|\nn P_\cdot f|$ is bounded on $[0,T]\times M$ for any $T>0$, then drop this assumption by an approximation argument.
 Since this condition is trivial for compact $M$, our proof below is much
 shorter than that in \cite{Hsu}.

 \beg{lem} \label{L2.1} Assume that $f\in C_b^1(M)$ such that $|\nn P_\cdot f|$ is bounded on
 $[0,T]\times M$ for any $T>0$. Then $(\ref{Hsu})$ holds. \end{lem}

 \beg{proof}  For any $\vv>0$, let

 $$\zeta_s=\ss{\vv +|\nn P_{t-s}f|^2}\, (X_s),\ \ \ s\le t.$$ By the
 It\^o formula we have

\beg{equation*}\beg{split} \d \zeta_s = & \d M_s+\ff{L|\nn
P_{t-s}f|^2 -2 \<\nn LP_{t-s}f, \nn P_{t-s}f\>}
{2\ss{\vv+|\nn P_{t-s}f|^2)^2}}(X_s) \d s\\
&- \ff{|\nn|\nn P_{t-s}f|^2|^2}{4(\vv+|\nn
P_{t-s}f|^2)^{3/2}}(X_s)\d s+\ff{N |\nn P_{t-s}f|^2}{2\ss{\vv+|\nn
P_{t-s}f|^2}} (X_s)\d l_s, \ \ \ s\le t,
\end{split}\end{equation*} where $M_s$ is a local martingale.  Combining this with (\ref{2.0}) and
(see \cite[(1.14)]{Ledoux})
\beq\label{Ledoux}L|\nn u|^2-2\<\nn Lu, \nn u\>\ge -2 \kk_1|\nn u|^2
+\ff{|\nn|\nn u|^2|^2}{2|\nn u|^2},\end{equation}   we obtain

$$\d\zeta_s\ge \d M_s -\ff{\kk_1 |\nn
P_{t-s}f|^2}{\vv+|\nn P_{t-s} f|^2}(X_s)\zeta_s\d s -\ff{\kk_2  |\nn
P_{t-s}f|^2}{\vv + |\nn P_{t-s}f|^2}(X_s)\zeta_s\d l_s,\ \ \ s\le
t.$$ Since   $\zeta_s$ is bounded on $[0,t],$ $\kk_1$ and $\kk_2$
are bounded, and by Lemma \ref{L2.2} $\EE \e^{\ll l_t}<\infty$ for
all $\ll>0$, this implies that

$$[0,t]\ni s\mapsto \zeta_s\exp\bigg[\int_0^s\ff{\kk_1|\nn P_{t-r}f|^2}{\vv+|\nn P_{t-r}f|^2}(X_r)\d r +\int_0^s
\ff{\kk_2|\nn P_{t-r}f|^2}{\vv+ |\nn P_{t-r} f|^2}(X_r)\d
l_r\bigg]$$ is a submartingale for any $\vv>0$. Letting
$\vv\downarrow 0$ we conclude that

$$[0,t]\ni s\mapsto |\nn P_{t-s}f|(X_s) \exp\bigg[\int_0^s\kk_1(X_r)\d r +\int_0^s\kk_2(X_r)\d
l_r\bigg]$$ is a submartingale as well. This completes the
proof.\end{proof}

 By Lemma \ref{L2.1}, to prove Proposition \ref{P2.1}  it suffices to confirm the boundedness
of $|\nn P_\cdot f|$ on $[0,T]\times M$ for $f\in C_b^1(M).$   Below
we first consider $f\in C_0^\infty(M)$ satisfying the Neumann
boundary condition.

\beg{lem}\label{L2.3} Assume {\bf (A)}. If $(\ref{C})$ holds then
for any $T>0$ and $f\in C_0^\infty(M)$ such that $Nf|_{\pp M}=0$,
$|\nn P_\cdot f|$ is bounded on $[0,T]\times M.$\end{lem}

\beg{proof} We shall take a conformal change of metric as in
\cite{W07} to make the boundary convex, so that the known  estimates
 for the convex case can be applied. As explained  on page 1436 in
 \cite{W07}, under assumption {\bf (A)} there exists $\phi\in
 C^\infty(M)$ and a constant $R>1$ such that $1\le\phi\le R, |\nn
 \phi|\le R, N\log\phi |_{\pp M}\ge\si,$ and $\nn\phi=0$ outside
 $\pp_r M.$ Since $\II\ge -\si$, by \cite[Lemma 2.1]{W07} $\pp M$ is
 convex under the new metric

 $$\<\cdot,\cdot\>=\phi^{-2}\<\cdot,\cdot\>.$$ Let $\DD', \nn',\Ric'$
 be corresponding to the new metric. By \cite[Lemma 2.2]{W07}

 $$L' :=\phi^2 L= \DD' + (d-2)\phi\nn\phi +\phi^2 Z=: \DD' +Z'.$$
 As in \cite{W08} we shall now calculate the curvature tensor $\Ric'-\nn' Z'$ under
 the new metric. By \cite[(9)]{W07}, for any unit vector $U\in TM$, $U':= \phi  U$
is unit under the new metric, and the corresponding Ricci curvature
satisfies

\beq\label{DD}\beg{split}\Ric'(U',U')\ge &\phi^2 \Ric(U,U)
+\phi\DD\phi -(d-3) |\nn
\phi|^2\\
&-2(U\phi)^2 + (d-2)\phi\Hess_\phi(U,U).\end{split}\end{equation}
 Noting that

$$\nn_X'Y= \nn_X Y- \<X, \nn \log \phi\>Y -\<Y,\nn \log\phi\>X +\<X,Y\>\nn\log \phi,\ \ \ X,Y\in TM,$$
we have

\beg{equation*} \beg{split}&\<\nn_{U'}Z', U'\>' = \<\nn_U Z', U\> -\<Z',\nn\log \phi\>\\
&=\phi^2 \<\nn_U Z,U\> +(U\phi^2) \<Z,U\> +(d-2)(U\phi)^2+(d-2)\phi
\Hess_{\phi} (U,U)-\<Z',\nn\log\phi\>.
\end{split}\end{equation*} Combining this with (\ref{DD}), (\ref{C}), $\|Z\|_r<\infty$
and the properties of $\phi$ mentioned above,
 we find a constant $K'\ge 0$ such that

$$\Ric'(U,'U')-\<\nn'_{U'}Z' ,U'\>'\ge -K',\ \ \ \<U',U'\>'=1.$$For
any $x,y\in M$, let $(X_t', Y_t')$ be the coupling by parallel
displacement of the reflecting diffusion processes generated by $L'$
with $(X_0',Y_0')=(x,y).$ Let $\rr'$ be the Riemannian distance
induced by $\<\cdot,\cdot\>'$. Since $(M, \<\cdot,\cdot\>')$ is
convex, we have (see \cite[(3.2)]{W97})

$$\rr'(X_t', Y_t')\le \e^{K't}\rr'(x,y),\ \ \ t\ge 0.$$ Since $1\le
\phi\le R$, we have $R^{-1}\rr\le \rr'\le\rr$ so that

\beq\label{L4} \rr(X_t',Y_t')\le R \e^{K't} \rr(x,y),\ \ \ t\ge
0.\end{equation} To derive the gradient estimate of $P_t$, we shall
make time changes

$$\xi_x(t)= \int_0^t \phi^2(X_s')\d s,\ \ \ \xi_y(t)=
\int_0^t\phi^2(Y_s')\d s.$$ Since $L'=\phi^2 L,$ we see that $X_t:=
X_{\xi_x^{-1}(t)}'$ and $Y_t:= Y'_{\xi_y^{-1}(t)}$ are generated by
$L$ with reflecting boundary. Again by $1\le \phi\le R$ we have

$$R^{-2} t\le \xi_x^{-1}(t), \xi_y^{-1}(t)\le t,\ \ \ t\ge 0.$$
Combining this with $|\nn \phi|\le R, 1\le \phi\le R$ and (\ref{L4})
we arrive at

\beg{equation}\label{L5}\beg{split}
|\xi_x^{-1}(t)-\xi_y^{-1}(t)|&\le \int_{\xi_x^{-1}(t)\land
\xi_y^{-1}(t)}^{\xi_x^{-1}(t)\lor \xi_y^{-1}(t)} \phi^2(Y_s')\d s=
|\xi_y\circ\xi_y^{-1}(t)-\xi_y\circ\xi_x^{-1}(t)|\\
&=|\xi_x\circ\xi_x^{-1}(t)-\xi_y\circ\xi_x^{-1}(t)|\le
\int_0^{\xi_x^{-1}(t)} |\phi^2(X_s')-\phi^2(Y_s')|\d s\\
&\le 2R^2 \rr(x,y) \int_0^t\e^{K's}\d s \le 2t\e^{K't} R^2
\rr(x,y).\end{split}\end{equation} Therefore,

\beq\label{L6} \beg{split} &|P_tf(x)-P_tf(y)| =
|\EE\{f(X_{\xi_x^{-1}(t)}')- f(Y'_{\xi_y^{-1}(t)})\}|\\
&\le \EE |f(X'_{\xi_y^{-1}(t)})-f(Y'_{\xi_y^{-1}(t)})| +
|\EE\{f(X'_{\xi_x^{-1}(t)})-f(X'_{\xi_y^{-1}(t)})\}|=:
I_1+I_2.\end{split}\end{equation} By (\ref{L4}) and $
\xi_y^{-1}(t)\le t$ we obtain

\beq\label{L7} I_1\le \|\nn f\|_\infty \e^{K' t}R
\rr(x,y).\end{equation} Moreover, since $f\in C_0^\infty(M)$ with
$Nf|_{\pp M}=0,$ it follows from the It\^o formula and (\ref{L5})
that

$$I_2\le \bigg|\EE\int_{\xi_x^{-1}(t)\land
\xi_y^{-1}(t)}^{\xi_x^{-1}(t)\lor \xi_y^{-1}(t)} L' f(X_s')\d
s\bigg|\le \|L'f\|_\infty \EE|\xi_x^{-1}(t)-\xi_y^{-1}(t)|\le
c_1t\e^{K't} \rr(x,y)$$ holds for some constant $c_1>0.$  Combining
this with (\ref{L6}) and (\ref{L7}) we conclude that

$$\|\nn P_tf\|_\infty\le c_2 (1+t) \e^{K't},\ \ \ t\ge 0$$ for some
constant $c_2>0$.\end{proof}

\ \newline\emph{Proof of Prposition \ref{P2.1}.}  Let $f\in C_b^1(M)$. By Lemma \ref{L2.1}
 we only have to prove the  boundedness
of $|\nn P_\cdot f|$ on $[0,T]\times M$ .

 (a) Let $f\in C_0^\infty(M).$ In this case there exist a sequence of functions
  $\{f_n\}_{n\ge 1}\subset C_0^\infty(M)$ such that $Nf_n|_{\pp M}=0, f_n\to f$ uniformly as $n\to\infty$,
  and $\|\nn f_n\|_\infty\le 1+\|\nn f\|_\infty$ holds
  for any $n\ge 1,$ see e.g. \cite{W94}. By Lemmas \ref{L2.1} and \ref{L2.3}, (\ref{Hsu})
  holds for $f_n$ in place of $f$ so that Lemma \ref{L2.2} implies

  $$\ff{|P_t f_n(x)-P_t f_n(y)|}{\rr(x,y)}\le C,\ \ \ t\le T, n\ge 1, x\ne y$$ for some constant $C>0.$
  Letting first $n\to 0$ then $y\to x$, we conclude that $|\nn P_\cdot f|$ is bounded
  on $[0,T]\times M.$

  (b) Let $f\in C_b^\infty(M)$.  Let $\{g_n\}_{n\ge 1}\subset C_0^\infty(M)$ be such that
  $0\le g_n\le 1, |\nn g_n|\le 2$ and $g_n\uparrow 1$ as $n\uparrow\infty$. By (a) and Lemma \ref{L2.1}, we may apply (\ref{Hsu}) to $g_nf$ in place of $f$
  such that Lemma  \ref{L2.2} implies

   $$\ff{|P_t (g_nf)(x)-P_t(g_n f)(y)|}{\rr(x,y)}\le C,\ \ \ t\le T, n\ge 1, x\ne y$$
   holds for some constant $C>0.$ By the same reason as in (a) we conclude that $|\nn P_\cdot f|$ is bounded
  on $[0,T]\times M.$

(c) Finally, for $f\in C_b^1(M)$ there exist $\{f_n\}_{n\ge
1}\subset C_b^\infty(M)$ such that $f_n\to f$ uniformly as
$n\to\infty$ and $\|\nn f_n\|_\infty\le \|\nn f\|_\infty+1$ for any
$n\ge 1.$ Therefore, the proof is complete by the same reason as in
(a) and (b). \qed

 \section{Proof of Theorem \ref{T1.1}}

 The boundedness of $P_t^Q$ under the uniform norm is ensured by Lemma
 \ref{L2.2}.   Since a bounded continuous function can be uniformly
approximated by bounded smooth functions, due to Lemma \ref{L2.2} we
may and do assume that $Q\in C_b^\infty(\pp M).$  To handle the integral $\int_0^t Q(X_s)\d l_s$,
we shall also need the upper bound of $L\rr_{\pp M}.$

\beg{lem}\label{L3.1} Let $\II\ge -\si$ and $(\ref{C})$ hold. Then

$$L\rr_{\pp M}\le (d-1)\si + \sup_{\pp M} \<Z, N\> + K \rr_{\pp M},\ \ \ \rr_{\pp M}<\i_{\pp M}.$$
\end{lem}

\beg{proof} Let $x\in M$ such that $\rr_{\pp M}(x)<\i_{\pp M}.$ Then there exist a unique $x_0\in\pp M$
and the minimal geodesic $x_\cdot: [0,\rr_{\pp M}(x)]\to M$ linking $\pp M$ and $x$. By (\ref{C})
we have

$$\Ric(\dot x_s, \dot x_s)\ge \<\nn_{\dot x_s}Z, \dot x_s\>-K=:R(s).$$
Let $h$ solve the equation

$$h''(s)+ \ff{R(s)}{d-1} h(s) =0,\ \ \ h(0)=1, h'(0)=\si.$$
By the Laplacian comparison theorem (see \cite[Theorem 1]{K2}),

\beg{equation*}\beg{split} \DD \rr_{\pp M}(x)&\le \ff{(d-1)h'}h (\rr_{\pp M}(x)) =
(d-1) \si +(d-1) \int_0^{\rr_{\pp M}(x)}\bigg(\ff{h'} h\bigg)'(s)\d s\\
&=(d-1)\si- (d-1) \int_0^{\rr_{\pp M}(x)} \Big(\ff{R(s)} {d-1} +\ff{(h')^2}{h^2}(s)\Big)\d s\\
&\le (d-1)\si +K\rr_{\pp M}(x)-\int_0^{\rr_{\pp M}(x)} \<\nn_{\dot x_s}Z, \dot x_s\>\d s.
\end{split}\end{equation*}
Then the proof is completed by noting that

$$Z\rr_{\pp M}(x)= \<Z,N\>(x_0) +\int_0^{\rr_{\pp M}(x)} \<\nn_{\dot x_s}Z, \dot x_s\>\d s.$$
\end{proof}

\subsection{The strong Feller property}

As explained after Theorem \ref{T1.1}, when $M$ is compact the
solution to (\ref{H}) is bounded in $x\in M$, so that (\ref{A0})
holds. In particular, $P_t^Qf$ is differentiable for $f\in \B_b(M)$
and thus, $P_t^Q$ is strong Feller. When $M$ is noncompact, this
argument does not apply due to the lack of boundedness of
$u(t,\cdot).$ Below we provide a different proof for the strong
Feller property.

 a) We first prove the Feller property. Since by
Lemma \ref{L2.2} for $f\in C_b(M)$ the function $P_t^Qf$ is bounded,
it suffices to show that

 \beq\label{C1} \lim_{y\to x} P_t^Qf(y)= P_t^Qf(x),\ \ \ x\in M.\end{equation}
 For any $y\in M$, let $(X_s, Y_s)$ be the coupling constructed in the proof of
 Lemma \ref{L2.3} via time changes. We shall first prove

 \beq\label{C2} \lim_{y\to x} \max_{s\in [0,t]} \rr(X_s, Y_s)=0.\end{equation}
 Using the notations in the proof of Lemma \ref{L2.3} and adopting (\ref{L4}) and (\ref{L5}),
there exists a constant $c(t)>0$ such that  for any $s\in [0,t]$

 \beg{equation*}\beg{split} &\rr(X_s, Y_s)= \rr(X'_{\xi_x^{-1}(s)}, Y'_{\xi_y^{-1}(s)})\le
\rr(X'_{\xi_y^{-1}(s)}, Y'_{\xi_y^{-1}(s)})
+\rr(X'_{\xi_x^{-1}(s)}, X'_{\xi_y^{-1}(s)})\\
&\le c(t) \rr(x,y) + \sup\big\{\rr(X_{s_1}', X_{s_2}'):\ s_1, s_2\in [0,t+
c(t)\rr(x,y)], |s_1-s_2|\le c(t)\rr(x,y)\big\}.\end{split}\end{equation*}
By the continuity of the reflecting diffusion process we prove (\ref{C2}).

 Next, to describe $\int_0^t Q(X_s) \d l_s$ we shall apply the It\^o
formula to a proper reference function of $X_s$.  To this end, we
first extend $Q$ to a smooth function on $M$. By assumption {\bf
(A)}, one may find a function  $\tt Q\in C^\infty(M)$ such that $\tt
Q|_{\pp M}=Q, N\tt Q|_{\pp M}=0$ and $|\nn \tt Q|+|L\tt Q|$ is
bounded. This can be realized by using the polar coordinates

$$\pp M\times [0,r)\ni (\theta, s)\mapsto \exp[sN_\theta]$$
for small enough $r>0$ such that $\rr_\pp$ is smooth on $\pp_r M$.
From this one may take  $\tt Q (\theta, s)= Q(\theta) h(s)$ on
$\pp_r M$ for some $h\in C^\infty([0,\infty)$ such that $h(0)=1,
h'(0)=0$ and $h(s)=0$ for $s\ge r,$ and let $\tt Q=0$ outside $\pp_r
M$. This $\tt Q$ meets our requirements since $L\rr_\pp$ is bounded on $\pp_r M$
according to (\ref{L0}) and  Lemma \ref{L3.1}.

 Let $\Phi\in
C_0^\infty([0,\infty))$ be such that $0\le\Phi\le 1, \Phi(s)=1$ for
$s\in [0,1]$ and $\Phi(s)=0$ for $s\ge 2.$ Let

$$\psi_n= \ff 1 n \int_0^{\rr_{\pp M}/n} \Phi(s)\d s.$$ Then $0\le
\psi_n\le 2n^{-1}, \psi_n=\rr_{\pp M}$ for $\rr_{\pp M} \le n^{-1},
\psi_n$ is constant for $\rr_{\pp M}\ge 2 n^{-1}$ and
$|\nn\psi_n|\le 1.$ Moreover, $\psi_n\in C^\infty(M)$ for large $n$.
Since $\nn \psi_n =N$ and $N\tt Q=0$ on $\pp M$, by the It\^o
formula we have

$$(\tt Q \psi_n)(X_t)= M_n(t) +\int_0^t L(\psi_n \tt Q)(X_s)\d s
+\int_0^t Q(X_s)\d l_s,$$ where $M_n(t)$ is a martingale with

\beq\label{C01}\<M_n\>(t)=\int_0^t |\nn (\tt Q \psi_n)|^2(X_s)\d
s.\end{equation}  Note that $L(\tt Q\psi_n)$ is bounded
since so is $|\nn \tt Q|+|L\tt Q|+ 1_{\pp_r M} |L\rr_{\pp M}|.$
Similarly, let $l_s^y$ be the local time of $Y_s$
on $\pp M$, we have

$$(\tt Q \psi_n)(Y_t)= M_n^y(t) +\int_0^t L(\psi_n \tt Q)(Y_s)\d s
+\int_0^t Q(Y_s)\d l_s^y$$ for some martingale $M_n^y(t)$  with

\beq\label{C02}\<M_n^y\>(t)=\int_0^t |\nn (\tt Q \psi_n)|^2(Y_s)\d
s.\end{equation} Combining these with (\ref{C2}) and using the
dominated convergence theorem, we obtain

\beq\label{C3} P_t^Qf(x) =\lim_{y\to x} \EE \big\{f(Y_t) \e^{(\tt Q
\psi_n)(Y_t) -\int_0^t L(\tt Q\psi_n)(Y_s)\d
s-M_n(t)}\big\}\end{equation} and

\beg{equation}\label{C4}\beg{split} I_n(y)&:= \big|P_t^Qf(y)- \EE
\big\{f(Y_t) \e^{(\tt Q \psi_n)(Y_t) -\int_0^t
L(\tt Q\psi_n)(Y_s)\d s-M_n(t)}\big\}\big|\\
&\le \EE\big|f(Y_t) \e^{\int_0^t Q(Y_s)\d l_s^y}(1-\e^{M_n^y(t)-M_n(t))}\big|\\
&\le (P_t^{2Q}f^2(y))^{1/2} (\EE|1-\e^{(M_n^y(t)-M_n(t))}|^2)^{1/2}\\
&\le
c_1(\EE|1-\e^{(M_n^y(t)-M_n(t))}|^2)^{1/2}\end{split}\end{equation}
for some constant $c_1>0.$   Since by the construction of $\psi_n$,
(\ref{C01}) and (\ref{C02}) we conclude that $M_n(t)$ and $M^y_n(t)$
are exponentially integrable uniformly in $y$ and

$$\<M_n\>(t)+\<M_n^y\>(t)\le c_2\int_0^t \{1_{\pp_{2/n}M}(X_s)+1_{\pp_{2/n}M}(Y_s)
\}\d s$$ holds for some constant $c_2>0$, it follows from (\ref{C2}) and (\ref{C4}) that

$$\lim_{n\to \infty}\lim_{y\to x} I_n(y)=0.$$
Combining this with (\ref{C3}) we derive (\ref{C1}).

b) Let $f\in \B_b(M).$ By Remark A.1 in the Appendix, the Neumann
semigroup is strong Feller. So, $f_\vv:= P_{\vv} f\in C_b(M)$ for
any $\vv>0.$ Combining this with the Feller property of $P_t^Q$,  it
suffices to prove

\beq\label{DD1} \lim_{\vv\to 0}
\|P_t^Qf-P_{t-\vv}^Qf_\vv\|_\infty=0.\end{equation} Since $Q$ is
bounded and $l_s$ is continuous in $s$ according to  (\ref{CC}) and
the continuity of $X_s$, Lemma \ref{L2.2} implies that

\beq\label{DD2}\lim_{\vv\to 0} \sup_{x\in
M}\EE^x\Big|f(X_t)\e^{\int_0^tQ(X_s)\d
l_s}\big(1-\e^{\int_{t-\vv}^tQ(X_s)\d
l_s}\big)\Big|=0.\end{equation} Next, let $\{\F_s\}_{s\ge 0}$ be the
natural filtration of $X_s$. By the Markov property we have

\beg{equation*} \beg{split} P_{t-\vv}^Qf_\vv &=\EE\big\{(P_\vv
f)(X_{t-\vv}) \e^{\int_0^{t-\vv} Q(X_s)\d l_s}\big\} \\
&= \EE \big\{\e^{\int_0^{t-\vv} Q(X_s)\d l_s} \EE
(f(X_t)|\F_{t-\vv})\big\}\\
&= \EE\big\{f(X_t) \e^{\int_0^{t-\vv} Q(X_s)\d
l_s}\big\}.\end{split}\end{equation*} Combining this with
(\ref{DD2}) we prove (\ref{DD1}).

 \subsection{The  generator}
We first prove that

\beq\label{D1} \lim_{t\to 0} \|P_t^Qf-f\|_\infty=0, \ \ \ f\in
C_0(M).\end{equation} Since a function in $C_0(M)$ can be uniformly
approximated by functions in $C_0^\infty(M)$ satisfying the Neumann
boundary condition, we may assume that $f\in C_0^\infty(M)$ with
$Nf|_{\pp M}=0.$ By the It\^o formula we have

$$\|P_tf-f\|_\infty \le \sup_M\int_0^t P_s|Lf|\d s\le t\|Lf\|_\infty$$ which goes to zero as $t\to 0$.
Noting that

$$|P_t^Qf-f|\le |P_t f-f|+ |P_t^Qf-P_tf|\le |P_tf-f|+ \|f\|_\infty \sup_{x\in M} \EE^x(\e^{\|Q\|_\infty l_t}-1),$$
by Lemma \ref{L2.2} for $\ll=\|Q\|_\infty$
 and letting first $t\to 0$ then $r\to 0$, we obtain (\ref{D1}).

Next, let $f\in C_0^2(M)$ satisfy the boundary condition
 $Nf+Qf=0.$ We intend to prove that

\beq\label{D}\lim_{t\to 0} \bigg\|\ff{P_t^Qf-f}t-Lf\bigg\|_\infty
=0.\end{equation}   By the It\^o
 formula and the boundary condition, we have

 \beg{equation*}\beg{split} &\d\big\{f(X_t) \e^{\int_0^t Q(X_s)\d
 l_s}\big\}\\
&= d M_t +(Lf(X_t))\e^{\int_0^t Q(X_s)\d l_s}\d t +\big\{Nf(X_t) +
(Qf)(X_t)\big\}\e^{\int_0^t Q(X_s)\d l_s} \d l_t\\
&= \d M_t +(Lf(X_t))\e^{\int_0^t Q(X_s)\d l_s}\d
t\end{split}\end{equation*} for some martingale $M_t$. This implies

$$P_t^Qf(x)= f(x)+ \int_0^t \EE^x\big\{(Lf(X_s))\e^{\int_0^s Q(X_r)\d l_r}\big\}
\d s.$$  So,
$$\bigg\|\ff {P_t^Qf-f} t-Lf\bigg\|_\infty \le \ff 1 t \int_0^t \|P_s^Q (Lf)- Lf\|_\infty\d s,$$
 and hence (\ref{D}) follows from (\ref{D1}) since $Lf\in C_0(M)$.

\subsection{The symmetry and $C_0$ property}  Let $Z=\nn V$.
Since  by Lemma \ref{L2.2}

$$(P_t^Q f)^2 \le (P_tf^2) \EE \e^{2l_t \|Q\|_\infty}\le c(t) P_t f^2$$ holds for some constant
$c(t)>0$, the boundedness of $P_t^Q$ in $L^2(\mu)$ follows from the
fact that $\mu$ is $P_t$-invariant. Moreover, since

\beg{equation*}\beg{split}|P_t^Qf-f|^2 &\le 2 |P_t^Q f- P_t f|^2 + 2|P_t f-f|^2\\
&\le 2|P_tf-f|^2 + 2(P_tf^2) \sup_{x\in M} \EE^x(\e^{\|Q\|_\infty
l_t}-1)^2,\end{split}\end{equation*} by Lemma \ref{L2.2} and the
strong continuity of $P_t$ in $L^2(\mu)$, we conclude that $P_t^Q$
is strongly continuous in $L^2(\mu)$ as well. So, it remains to
prove that for any $f,g\in C_0(M)$

\beq\label{N1} \mu(gP_t^Qf)= \mu(fP_t^Qg).\end{equation}
 We shall prove (\ref{N1}) by using
symmetric Schr\"odinger semigroups to approximate $P_t^Q$.

Let $\tt Q$ and $\psi_n$ be constructed above.  We have
\beq\label{3.1} \EE^x [f(X_t)\e^{\int_0^tQ(X_s)\d l_s}] =
\EE^x[f(X_t)\e^{- \int_0^t L (\tt Q\psi_n)(X_s)\d s}] +
\vv_n,\end{equation} where

$$\vv_n:= \EE^x [f(X_t) \e^{-\int_0^tQ(X_s)\d l_s}(1-\e^{- M_n(t) -
(\tt Q\psi_n)(X_t)})],$$  which goes to zero uniformly in $x$ as
$n\to\infty$ according to Lemma \ref{L2.2} and the properties of
$\tt Q$ and $\psi_n$. Let $P_t^{(n)}$ be the Schr\"odinger semigroup
generated by

$$L_n:= L-  L(\tt Q\psi_n).$$  Since  $L(\tt Q\psi_n)$ is bounded,
by the Feynman-Kac formula

$$P_t^{(n)}f(x)= \EE^x[f(X_t)\e^{- \int_0^tL (\tt Q\psi_n)(X_s)\d s}]$$
and $P_t^{(n)}$ is symmetric in $L^2(\mu)$. So,

\beq\label{3.2} \beg{split}&\int_M
g(x)\EE^x[f(X_t)\e^{-\int_0^tL(\tt Q\psi_n)(X_s)\d s}]\mu(\d x)=
\int_M
g(x) P_t^{(n)} f(x)\mu(\d x)\\
&=\int_Mf(x)P_t^{(n)}g(x)\mu(\d x) = \int_M f(x)
\EE^x[g(X_t)\e^{-\int_0^tL(\tt Q\psi_n)(X_s)\d s}]\mu(\d x)\\
&=\int_M f(x)\EE^x[g(X_t)\e^{\int_0^tQ(X_s)\d l_s +M_n(t)-(\tt
Q\psi_n)(X_t)}]\mu(\d x).\end{split}\end{equation} Obviously,
$M_n(t)-(\tt Q\psi_n)(X_t)\to 0$ a.s. as $n\to \infty$ and is
exponentially integrable uniformly in $x$. So, by Lemma \ref{L2.2},
 (\ref{3.1}), (\ref{3.2}) and using the dominated convergence theorem we arrive at

$$\mu(g P_t^Q f) = \lim_{n\to\infty} \bigg\{\int_M f(x)
\EE^x[g(X_t)\e^{\int_0^t Q(X_s)\d l_s +M_n(t)-(\tt
Q\psi_n)(X_t)}]\mu(\d x)+\vv_n\bigg\}=\mu(fP_t^Qg).$$

\subsection{The Dirichlet form}

Again let $Z=\nn V$ for $V\in C^2(M)$. Let $Q\le 0$ such that $\E\ge
0.$ Since by (\ref{Robin}) and the integration by parts formula, we
have

$$\E(f,g)= -\int_M fL g\d\mu,\ \ \ f,g\in \D_0,$$ the form
$\E,\D_0)$ is closable. Moreover, as in \S3.2 for $f\in \D_0$ one
has $\ff{P_t^Qf-f}t \to Lf$ in $L^2(\mu)$ as $t\to 0$, it remains to
show that $(\E,\D_0)$ is  a pre-Dirichlet form in $L^2(\mu)$.
Firstly, to understand that $(\E,\D_0)$ is  well defined   in
$L^2(\mu)$, i.e. $\E(f,g)$ is independent of $\mu$-versions of $f$
and $g$, for a bounded continuous
  extension $\tt Q$ of $Q$ on $M$, we rewrite

$$\mu_\pp (Qfg)= \lim_{r\to 0} \ff {\mu(\tt Q f g 1_{\pp_r
M})}{r}.$$ Since $Q\le 0$ implies the nonnegativity  and the normal
contraction property of $\E$, it remains to show that $\D_0$ is
dense in $L^2(\mu).$

It is well known that the class of functions in $C_0^\infty(M)$
satisfying the Neumann boundary condition is dense in $L^2(\mu)$, it
suffices to prove that for any $f\in C_0^\infty(M)$ with $Nf|_{\pp
M}=0,$ there exists a sequence $\{f_n\}\subset \D_0$ such that
$\mu(|f_n-f|^2)\to 0$ as $n\to\infty.$ To this end, for any $\vv>0$,
let $h_\vv \in C_0^\infty([0,\infty))$ such that $h_\vv (0)=0,
h_\vv'(0)=1, h_\vv(s)=0$ for $s\ge r_0$ and $||h_\vv||_\infty\le
\vv/\|Q\|_\infty.$ Here $r_0>0$ is such that $\rr_{\pp M}$ is smooth
on $\pp_{r_0}M.$  Then $\psi_\vv (\theta, s):= 1+ Q(\theta)h_\vv(s)$
defined under the polar coordinates

$$ \pp M \times [0,r_0) \ni
(\theta, s)\mapsto \pp_{r_0}M$$ is smooth and can be naturally
extended smoothly on $M$ by letting $\psi_\vv=1$ on $M\setminus
\pp_{r_0}M.$ Obviously, we have $\psi_\vv|_{\pp M}=1,
N\psi_\vv|_{\pp M}=Q$ and $|\psi_\vv-1|\le \vv.$  Thus, $\psi_\vv
f\in \D_0$ and $\psi_\vv f\to f$ in $L^2(\mu)$ as $\vv\to 0.$

\section{Proof of Theorem \ref{T1.2}}

 By Lemma \ref{L2.2},
it remains to verify (\ref{HWI2}). Let $f\in C_b^1(M)$ and $t>0$. We
have

\beq\label{3.3} \ff{\d}{\d s} P_s\big\{(P_{t-s}f^2)\log
P_{t-s}f^2\big\}=P_s\ff{|\nn P_{t-s}f^2|^2}{P_{t-s}f^2},\ \ \ s\in
[0,t].\end{equation} By (\ref{G2}) and the Schwartz inequality we
have

\beg{equation*}\beg{split} \ff{|\nn
P_{t-s}f^2|^2}{P_{t-s}f^2}(y)&\le \e^{2K(t-s)}\ff{(\EE^y\{|\nn
f^2|(X_{t-s})\e^{\si l_{t-s}}\})^2}{P_{t-s} f^2(y)}\\
&\le 4 \e^{2K(t-s)} \EE^y \{|\nn f|^2(X_{t-s})\e^{2\si
l_{t-s}}\}=:4\e^{2K(t-s)}g_s(y),\ \ \ s\in [0,t], y\in
M.\end{split}\end{equation*} Combining this with (\ref{3.3}) we
obtain

$$P_t (f^2\log f^2)\le (P_tf^2)\log P_t f^2 + 4 \int_0^t
\e^{2K(t-s)}P_s g_s\d s.$$Since $\mu$ is an invariant measure of
$P_t$, taking integral for both sides with respect to $\mu$ we
arrive at

$$\mu(f^2\log f^2) \le \mu((P_tf^2)\log P_t f^2)+ 4\int_0^t
\e^{2K(t-s)}\mu(g_s)\d s.$$ Since by Theorem \ref{T1.1}(2)
$P_{t-s}^Q$ for $Q= 2\si$ is symmetric in $L^2(\mu)$, we have

$$\mu(g_s)= \int_M |\nn f|^2(x)\EE^x\e^{2\si l_{t-s}}\mu(\d x)\le
\mu(|\nn f|^2) \eta_{2\si}(t-s),$$ it follows that

\beq\label{3.4} \mu(f^2\log f^2) \le \mu((P_tf^2)\log P_t f^2)+
4\mu(|\nn f|^2 ) \int_0^t \e^{2Ks}\eta_{2\si}(s)\d s.\end{equation}
This is an extension of (\ref{LS2}) to the nonconvex case.

On the other hand, we intend to establish an analogous to (\ref{W})
for the present situation. For any $x,y\in M$, let $x_\cdot:
[0,1]\to M$ be the minimal curve linking $x$ and $y$ with constant
speed. We have $|\dot x_s|=\rr(x,y)$. Let $h\in C^1([0,t])$ be such
that $h_0=1, h_t=0.$ Then by (\ref{G2}) which follows from
Proposition \ref{P2.1}, we have

\beg{equation}\label{3.5}\beg{split} &P_t\log f^2(x)-\log P_t
f^2(y)= \int_0^t
\ff{\d }{\d s} P_s(\log P_{t-s}f^2)(x_{h_{t-s}})\d s\\
&\le \int_0^t \Big\{|\dot h_{t-s}| \rr(x,y) |\nn P_s(\log
P_{t-s}f^2)|(x_{h_{t-s}})-\EE^{x_{h_{t-s}}} \ff{|\nn
P_{t-s}f^2|^2}{(P_{t-s}f^2)^2}(X_s)  \Big\}\d s\\
&\le \int_0^t \EE^{x_{h_{t-s}}} \Big\{ |\dot h_{t-s}| \rr(x,y)
\ff{|\nn P_{t-s}f^2|}{P_{t-s}f^2}(X_s)\e^{K(t-s)+\si
l_{t-s}}-\ff{|\nn
P_{t-s}f^2|^2}{(P_{t-s}f^2)^2}(X_s)\Big\}\d s\\
&\le \ff {\rr(x,y)^2}4 \int_0^t {\dot h}_s^2
\e^{2Ks}\eta_{2\si}(s)\d s=:c(t)\rr(x,y)^2.\end{split}\end{equation}

Now, let $\mu(f^2)=1$ and $\pi\in \scr C(f^2\mu,\mu)$ be the optimal
coupling for $W_2(f^2\mu,\mu).$ It follows from the symmetry of
$P_t$ and (\ref{3.5}) that

\beg{equation*}\beg{split} \mu((P_tf^2)\log P_tf^2)&=\mu(f^2 P_t\log
P_t f^2)=\int_{M\times M} P_t(\log P_t f^2)(x)\pi(\d x,\d y)\\
&\le \int_{M\times M} \big\{\log
P_{2t}f^2(y)+c(t)\rr(x,y)^2\big\}\pi(\d x,\d y)\\
&=\mu(\log P_{2t}f^2)+c(t)W_2(f^2\mu,\mu)^2\le
c(t)W_2(f^2\mu,\mu)^2,\end{split}\end{equation*} where in the last
step we have used the Jensen inequality that

$$\mu(\log P_{2t}f^2)\le \log \mu(P_{2t} f^2)=0.$$ Combining this with
(\ref{3.4}) we obtain

$$\mu(f^2\log f^2) \le 4\mu(|\nn f|^2) \int_0^t
\e^{2Ks}\eta_{2\si}(s)\d s+ \ff{W_2(f^2\mu,\mu)^2} 4 \int_0^t {\dot
h}_s^2 \e^{2Ks} \eta_{2\si}(s)\d s.$$ Then the  proof is completed
by taking

$$h_s= \ff{\int_s^t \e^{-2Ku}\eta_{2\si}(u)^{-1}\d u}{\int_0^t
\e^{-2Ku}\eta_{2\si}(u)^{-1}\d u},\ \ \ s\in [0,t].$$ \qed

\section{Appendix: the Bismut formula}

By using a formula for the gradient of $P_t$ derived in \cite{Hsu},
one obtains the following Bismut type formula (\ref{B}) as in
\cite{T}, which in particular implies the strong Feller property of
$P_t$ as explained in Remark A.1 below.

Because of the exponential integrability of $l_t$ ensured by Lemma
\ref{L2.2}, it is easy to see that the argument in \cite{Hsu} for
compact $M$ works also for the present case under assumption {\bf
(A)} and condition (\ref{C}). To state the  formula for the gradient
of $P_t$ obtained in \cite{Hsu}, let us first introduce the SDE for
the horizontal lift of the reflecting $L$-diffusion process.

Let $O(M)$ be the bundle of orthonormal frames over $M$ and let
$\pi: O(M)\to M$ be the natural projection. Then $X_t$ and its
horizontal lift $u_t$ on $O(M)$ solve the following equations:

\beg{equation*}\beg{split} & \d u_t= H_{u_t} \circ\d X_t,\\
&\d X_t= \ss 2\, u_t\circ \d B_t + Z(X_t)\d t+ N(X_t)\d
l_t,\end{split}\end{equation*}where $B_t$ is the $d$-dimensional
Brownian motion and $H_u: T_{\pi u}M\to T_u O(M)$ is the horizontal
lift at $u\in O(M).$ Next, let $\mathbb M_t$ be the
$\R^d\otimes\R^d$-valued process solving the equation
$$ \d \mathbb M_t= -\mathbb M_t R_{u_t}\d t,\ \ \ \mathbb M_0=I,$$ where

$$R_u(a,b)= \Ric(ua, ub)-\<\nn_{ua}Z, ub\>,\ \ \ u\in O(M), a,b\in
\R^d.$$ Then for any $f\in C_0^\infty(M)$ we have (cf. the proof of
\cite[Theorem 5.1]{Hsu})

\beq\label{App1} u_0^{-1} \nn P_t f= \EE \big\{\mathbb
M_tu_t^{-1}\nn f(X_t)\big\}.\end{equation} Since by Lemma \ref{L2.2}
and Proposition \ref{P2.1} $|\nn P_\cdot f|$ is bounded on
$[0,T]\times M$ for any $T>0,$ this follows since according to
\cite[Theorem 3.7]{Hsu} the process  $\{M_s F(u_s, t-s)\}_{s\in
[0,t]}$ is a martingale for $F(u,s)= u^{-1} \nn P_s f(\pi u).$ Due
to (\ref{App1}) we have the following result on the Bismut type
formula.

\beg{prp} Assume {\bf (A)} and $(\ref{C})$. Then for any $f\in
C_0^\infty(M)$ and any increasing function $h\in C^1([0,t])$ such
that $h(0)=0, h(t)=1,$

\beq\label{B} \nn P_t f(x)= \ff{u_0}{\ss 2} \EE \bigg[f(X_t)
\int_0^t h'(s) \<\mathbb M_s, \d B_s\>\bigg]\end{equation}holds for
$x\in M$ and $X_t, u_t$ start from $x, u_0\in O_x(M)$ respectively.
\end{prp}

\beg{proof} By It\^o's formula we have

$$\d P_{t-s} f(X_s) =\ss 2\, \<\nn P_{t-s}f(X_s), u_s\d B_s\>.$$
Then

\beq\label{App3} f(X_t)= P_t f+\ss 2\, \int_0^t\<\nn P_{t-s}f(X_s),
u_s \d B_s\>.\end{equation} Combining this with (\ref{App1}), for
any $a\in \R^d$,

\beg{equation*}\beg{split} &\ff 1 {\ss 2} \EE \bigg\{f(X_t)\int_0^t
h'(s) \<\mathbb M_s a, \d B_s\> \bigg\}= \EE \int_0^t \<\nn
P_{t-s}f(X_s), \mathbb M_s
a\>h'(s)\d s\\
&= \int_0^t \<u_0a, \nn P_t f\>h'(s)\d s= \<u_0a, \nn P_t
f\>.\end{split}\end{equation*}This completes the proof since $a\in
\R^d$ is arbitrary. \end{proof}

\paragraph{Remark A.1.} By (\ref{C}) and letting $\II\ge -\si,$  we
have

$$\|M_s\|\le \e^{Ks+\si l_s},\ \ \ s\ge 0.$$ So, by Lemma \ref{L2.2} and (\ref{B}), for any $t>0$ there exists a constant $C(t)>0$ such that

$$\|\nn P_t f\|_\infty\le \|f\|_\infty C(t),\ \ \ t>0, f\in
C_0^\infty(M).$$ This implies

$$|P_t f(x)-P_t f(y)|\le C(t) \|f\|_\infty \rr(x,y),\ \ \ x,y\in M,
f\in C_0^\infty(M).$$ By the monotone class theorem, this inequality
holds indeed for all $f\in \B_b(M)$ and thus, $P_t$ is  strong
Feller.

  \beg{thebibliography}{99}

 \bibitem{BE}   D. Bakry and M. Emery, \emph{Hypercontractivit\'e de
semi-groupes de diffusion}, C. R. Acad. Sci. Paris. S\'er. I Math.
299(1984), 775--778.

\bibitem{BGL}  S. G. Bobkov, I.  Gentil and M.  Ledoux,
\emph{Hypercontractivity of Hamilton-Jacobi equations,} J. Math.
Pures Appl. 80(2001), 669--696.

\bibitem{OV} F.  Otto and C. Villani, \emph{Generalization of an inequality
by Talagrand and links with the logarithmic Sobolev inequality,} J.
Funct. Anal.  173(2000), 361--400.

\bibitem{O} E. M. Ouhabaz, \emph{Analysis of Heat Equations on
Domians,} Princeton University Press, Princeton, 2005.

\bibitem{Hsu} E. P. Hsu, \emph{Multiplicative functional for the heat
equation on manifolds with boundary,} Michigan Math. J. 50(2002),
351--367.

\bibitem{K2} A. Kasue, \emph{On Laplacian and Hessian comparison theorems,}
Proc. Japan Acad. 58(1982), 25--28.

 \bibitem{Kasue} A. Kasue, \emph{Applications of Laplacian and Hessian comparison theorems,} Geometry of
 geodesics and related topics (Tokyo, 1982), 333--386, Adv. Stud. Pure Math., 3, North-Holland, Amsterdam, 1984.

\bibitem{Ledoux} M. Ledoux, \emph{The geometry of Markov diffusion generators,} Ann. Facu. Sci. Toulouse
9(2000), 305--366.

\bibitem{Q} Z. Qian, \emph{A gradient estimate on a manifold with convex boundary,}
 Proc. Roy. Soc. Edinburgh Sect. A 127(1997), 171--179.

 \bibitem{T} A. Thalmaier, \emph{On the differentiation of the heat
 semigroups and Poisson integrals,} Stochastics 61(1997), 297--321.

 \bibitem{W94} F.-Y. Wang, \emph{Application of coupling methods to the Neumann eigenvalue
  problem}, Probab. Theory Related Fields 98 (1994),  299--306.

\bibitem{W97} F.-Y. Wang, \emph{On estimation of the logarithmic Sobolev constant and
gradient estimates of heat
semigroups,} Probab. Theory Relat. Fields 108(1997), 87--101.

 \bibitem{W05} F.-Y. Wang, \emph{Gradient estimates and the first
Neumann eigenvalue on manifolds with boundary,} Stoch. Proc. Appl.
115(2005), 1475--1486.

 \bibitem{W07} F.-Y. Wang, \emph{Estimates of the first Neumann eigenvalue and the log-Sobolev constant
  on  nonconvex manifolds,} Math. Nachr. 280(2007), 1431--1439.

 \bibitem{W08} F.-Y. Wang, \emph{Gradient and Harnack inequalities on noncompact manifolds
 with boundary,}  preprint.

\end{thebibliography}

\end{document}